\documentclass[leqno]{amsart}
\usepackage[normal]{boilerart}
\usepackage{mathrsfs}
\usepackage{tikz-cd}

\swapnumbers
\newtheorem{mainthm}{Theorem}

\newcommand{\cop}{\copyright}

\newcommand{\Proj}{\categ{Proj}}

\newcommand{\Vect}{\categ{Vect}}

\def\reflectedop#1#2{\mathop{\reflectbox{$#1{#2}$}}}

\title[A general introduction]{Parametrized higher category theory and higher algebra: {A} general introduction}

\author{Clark Barwick}
\address{Department of Mathematics, Massachusetts Institute of Technology, 77 Massachusetts Avenue, Cambridge, MA 02139-4307, USA}
\email{clarkbar@math.mit.edu}

\author{Emanuele Dotto}
\address{Department of Mathematics, Massachusetts Institute of Technology, 77 Massachusetts Avenue, Cambridge, MA 02139-4307, USA}
\email{dotto@math.mit.edu}

\author{Saul Glasman}
\address{University of Minnesota, School of Mathematics, Vincent Hall, 206 Church St. SE, Minneapolis, MN 55455, USA}
\email{saulglasman0@gmail.com}

\author{Denis Nardin}
\address{Department of Mathematics, Massachusetts Institute of Technology, 77 Massachusetts Avenue, Cambridge, MA 02139-4307, USA}
\email{nardin@math.mit.edu}

\author{Jay Shah}
\address{Department of Mathematics, Massachusetts Institute of Technology, 77 Massachusetts Avenue, Cambridge, MA 02139-4307, USA}
\email{jshah@math.mit.edu}

\begin{document}

\maketitle

Let $k$ denote a field, and let $E\supseteq k$ be a finite Galois extension thereof with Galois group $G$. The algebraic $K$-groups $K_n(k)$ and $K_n(E)$, as defined by Quillen, together exhibit some interesting structure. Since these groups are defined in terms of the categories of finite-dimensional vector spaces (along with their additive structure), the forgetful functor $\fromto{\Vect(E)}{\Vect(k)}$ and the functor $\fromto{\Vect(k)}{\Vect(E)}$ given by $\goesto{X}{X\otimes_kE}$ give rise to homomorphisms
\[V\colon\fromto{K_n(E)}{K_n(k)}\text{\quad and\quad}F\colon\fromto{K_n(k)}{K_n(E)}.\]
Ordinary Galois theory shows that the composite functor $\fromto{\Vect(E)}{\Vect(E)}$ given by $\goesto{Y}{Y\otimes_kE}$ can be described as the direct sum
\[\bigoplus_{g\in G}g\colon\fromto{\Vect(E)}{\Vect(E)},\]
where $G$ acts in the obvious manner. Accordingly, we have an action of $G$ on $K_n(E)$ for which both $V$ and $F$ are equivariant, and a formula
\[FV=\sum_{g\in G}g.\]

Note that the equivariance of $V$ implies that it factors through the orbits $K_n(E)_G$, and the equivariance of $F$ implies that it factors through the fixed points $K_n(E)^G$, but these maps do not typically identify $K_n(k)$ with either the orbits or the fixed points. The data of $K_n(k)$ is an added piece of structure; that is, $K_n(k)$ cannot in general be recovered from $K_n(E)$ as a $G$-module.

But the problem is even deeper than this. Even if one considers all the $K$-groups together as a single entity (by thinking of these groups as the homotopy groups of a space or spectrum), one can construct a descent spectral sequence
\[E^2_{p,q}=H^{-p}(G,K_q(E)),\]
but this will not, as a rule, converge to the groups $K_{p+q}(k)$. In other words, the space or spectrum $K(k)$ is not the homotopy fixed point space/spectrum of the action of $G$ on the space/spectrum $K(E)$. Consequently, even knowing the homotopy type $K(E)$ with its action of $G$ is insufficient to recover the groups $K_n(k)$. This is the \emph{descent problem in algebraic $K$-theory}.

There is, of course, no need to consider the $K$-theories of $E$ and $k$ in isolation. One can also include the information of the $K$-groups of all the various subextensions $E\supseteq L\supseteq k$. In other words, for any subgroup $H\leq G$, one can contemplate the $K$-groups $K_n(E^H)$ of the fixed field $E^H$. These abelian groups each have conjugation homomorphisms
\[c_g\colon\fromto{K_n(E^H)}{K_n(E^{gHg^{-1}})}\]
for any $g\in G$. Additionally, for subgroups $K,L\leq H\leq G$, one again has the forgetful functor $\fromto{\Vect(E^K)}{\Vect(E^H)}$ and the functor $\fromto{\Vect(E^H)}{\Vect(E^L)}$ given by $\goesto{Y}{Y\otimes_{E^H}E^L}$, so again one has homomorphisms
\[V_K^H\colon\fromto{K_n(E^K)}{K_n(E^H)}\text{\quad and\quad}F_L^H\colon\fromto{K_n(E^H)}{K_n(E^L)}.\]
Again, a small amount of Galois theory reveals that these two homomorphisms compose in the following manner:
\[F_L^HV_K^H=\sum_{x\in L\backslash G/K}V_{L\cap(xKx^{-1})}^Lc_xF_{(x^{-1}Lx)\cap K}^K\colon\fromto{K_n(E^K)}{K_n(E^L)}.\]
And again, of course, the groups $K_n(E^H)$ cannot be recovered from the $G$-module $K_n(E)$ or the homotopy type $K(E)$ with its action of $G$.

Combined, this structure on the assignment $\goesto{H}{K_n(E^H)}$ makes up what is called a \emph{Mackey functor} for $G$. As we see, this is strictly more structure than a $G$-module. Similarly, the assignment $\goesto{H}{K(E^H)}$ is a \emph{spectral Mackey functor} for $G$ in the sense of the first author \cite{M1}. This is strictly more structure than a spectrum with a $G$-action. We call this object the \emph{$G$-equivariant $K$-theory of $E$ over $k$}.

In this monograph, we tease out the kind of structure on the categories $\Vect(E^H)$ that provides their $K$-theory with the structure of a spectral Mackey functor for $G$. As a first approximation, we note that, because the category of subextensions of $E$ is equivalent to the category of transitive $G$-sets, the functors $\goesto{Y}{Y\otimes_{E^H}E^L}$ together define what we call a \emph{$G$-category} -- a diagram of categories indexed on the opposite of the orbit category $\OO_G$ of $G$. Let us write $\underline{\Vect}_{E\supseteq k}$ for this $G$-category.

Of course, the $G$-category $\underline{\Vect}_{E\supseteq k}$ is relatively simple: after all, if one thinks of the action of $G$ on $\Vect(E)$, then $\Vect(E^H)$ is the category of $E$-vector spaces equipped with a semilinear action of $H$. In other words, $\Vect(E^H)$ is simply the homotopy fixed point category for the action of $H$ on $\Vect(E)$. So we might at first contemplate $\Vect(E)$ with its $G$-action. However, the adjoints to the functors in this $G$-category -- the forgetful functors -- contain extra information that compels us to contemplate entire $G$-category structure.

For example, the forgetful functor $\fromto{\Vect(E)}{\Vect(k)}$ is a kind of generalized \emph{product} of vector spaces: we regard it as \emph{indexed}, not over a mere set, but over the $G$-set $G/e$. To see why this is appropriate, first note that by the normal basis theorem, if $Y$ is an $E$-vector space with basis $\{v_i\}_{1\leq i\leq n}$, then there is an element $\theta\in E$ such that $Y$ has basis $\{g\theta v_i\}_{1\leq i\leq n,g\in G}$ over $k$. But without choosing this element, we would still be entitled to write
\[\prod_{\alpha\in G/e}Y\]
for this $k$-vector space. In the same manner, the presence of all the other right adjoints $\fromto{\Vect(E^H)}{\Vect(E^K)}$ in this diagram of categories can be regarded as the existence of various \emph{indexed products}
\[\prod_{\alpha\in K/H}Z\]
on this $G$-category. At the same time, since our field extensions are separable, these right adjoints are all also \emph{left} adjoints, and so we even think of this as endowing our $G$-category with \emph{indexed direct sums}
\[\bigoplus_{\alpha\in K/H}Z.\]

The point here is that the transfer structure on the equivariant algebraic $K$-groups arises from the additional structure of indexed products or coproducts on the $G$-category $\underline{\Vect}_{E\supseteq k}$. And this example refects a general principle: to get the full structure of a Mackey functor on equivariant algebraic $K$-theory of $E$ over $k$, one must work not only with the diagram of categories indexed by $\OO_G^{\op}$, but also the $G$-direct sums thereupon.

The $G$-category $\underline{\Vect}_{E\supseteq k}$ also carries a sophisticated multiplicative structure. Of course, the tensor product over $k$ provides an \emph{external product}
\[
\fromto{\Vect(E^K)\times\Vect(E^L)}{\Proj^{\,\textit{fg}}(E^K\otimes_k E^L)\simeq\prod_{x\in L\backslash G/K}\Vect(E^{(x^{-1}Lx)\cap K})}.
\]
In \cite{M2}, we demonstrated that the external products provide the equivariant algebraic $K$-groups with the structure of a graded Green functor, and, even better, they provide the equivariant algebraic $K$-theory spectra with the structure of a spectral Green functor.

However, there is a still richer multiplicative structure, whose impact on equivariant $K$-theory is studied here for the first time. Just as the usual norm of an element of $E$ is automatically Galois-invariant, we see that for any finite-dimensional $E$-vector space $V$, the tensor power $V^{\otimes G}$ comes with canonical descent data. We call the resulting $k$-vector space $N_E^k(V)$ the \emph{multiplicative norm} from $E$ to $k$. Quite simply, $N_E^k(V)$ is the $k$-vector space (of dimension $(\dim V)^{\# G}$) such that the set $\Hom_k(N_E^k(V),W)$ is in bijection with the set of \emph{norm forms $\fromto{V^{\times G}}{W\otimes_k E}$ for $E/k$} -- i.e., $k$-multilinear maps
\[\Phi\colon\fromto{V^{\times G}}{W\otimes_k E}\]
such that for any element $(v_h)_{h\in G}\in V^{\times G}$, any element $g\in G$, and any element $\lambda\in E$,
\[\Phi((v'_h)_{h\in G})=(g\lambda)\Phi((v_h)_{h\in G}),\]
where
\[
v'_h=\begin{cases}
\lambda v_g&\textup{if }h=g;\\
v_h&\textup{if }h\neq g,
\end{cases}
\]
and
\[
g\Phi((v_h)_{h\in G})=\Phi((v_{gh})_{h\in G}.
\]
So, $N_E^k(V)$ is the dual of the $k$-vector space of norm forms $\fromto{V^{\times G}}{E}$ for $E/k$. In particular, when $k=\RR$ and $E=\CC$, then $N_{\CC}^{\RR}(V)$ is precisely the dual space of the $\RR$-vector space of hermitian forms on $V$.

More generally, there are multiplicative norms for any subgroups $K\leq L\leq G$. Together with the external products, these multiplicative norms furnish $\underline{\Vect}_{E\supseteq k}$ with a \emph{$G$-symmetric monoidal structure}. In effect, this provides tensor products indexed over any finite $G$-set $U=\coprod_{i\in I}(G/H_i)$, which amount to functors
\[\bigotimes_{u\in U}\colon\fromto{\prod_{i\in I}\Vect(E^{H_i})}{\Vect(k)},\]
which are suitably associative and commutative.

This additional structure on $\underline{\Vect}_{E\supseteq k}$ descends to an analogous structure on the equivariant algebraic $K$-theory of $E$ over $k$. These provide the equivariant algebraic $K$-theory of $E$ over $k$ with the full structure of a \emph{$G$-$E_{\infty}$-algebra}.

\section*{Hill's program} To tell this story, we pursue here the general theory of \emph{$G$-$\infty$-categories}. But we are by no means the first to contemplate this possibility.

In their landmark solution of the Kervaire Invariant Problem \cite{MR3505179}, Mike Hill, Mike Hopkins, and Doug Ravenel developed a perspective on equivariant stable homotopy thery that centered on the study of indexed products, indexed coproducts, and indexed symmetric monoidal structures (incorporating their multiplicative norms). They argued that these structures were fundamental to the basic structure of equivariant stable homotopy theory.

In 2012, Hill presented (partly jointly with Hopkins) a sketch of a program to rewire huge swaths of higher category theory in order to embed these structures into the very fabric of the homotopy theory of homotopy theories. Hill sought a theory of $G$-$\infty$-categories and $G$-functors, along with a concomitant theory of internal homs, $G$-limits, $G$-colimits, $G$-Kan extensions, etc. He furthermore conjectured that, equipped with this technology, one could prove the following, which is an analogue of the universal property of the $\infty$-category $\Top$ of spaces.

\begin{mainthm}\label{TopGuniv} The $G$-$\infty$-category $\underline{\Top}_G$ of $G$-spaces -- whose value on an orbit $G/H$ is the $\infty$-category of $H$-spaces -- is freely generated under $G$-colimits by the contractible $G$-space; that is, for any $G$-$\infty$-category $D$ with all $G$-colimits, evaluation on the generator defines an equivalence of $G$-$\infty$-categories
\[\equivto{\underline{\Fun}_G^L\left(\underline{\Top}_G,D\right)}{D}.\]
Here $\underline{\Fun}_G^L$ is the $G$-$\infty$-category of $G$-colimit-preserving functors.
\end{mainthm}

\noindent In this text, we develop all this machinery, and this is the first main theorem.

Recall that one may speak of semiadditive $\infty$-categories, in which finite products and finite coproducts exist and coincide. In the same manner, Hill expected that one may speak of $G$-semiadditive $\infty$-categories, in which finite $G$-products and finite $G$-coproducts exist and coincide. Furthermore, the effective Burnside $\infty$-category $A^{\eff}(\FF)$ of finite sets is equivalent to the $\infty$-category of finitely generated free $E_{\infty}$-spaces, whence it is the free semiadditive $\infty$-category on one generator. Accordingly, in equivariant higher category theory, we have the following.

\begin{mainthm}\label{Aeffuniv} The $G$-$\infty$-category $\underline{A}^{\eff}(G)$ -- whose value on $G/H$ is the effective Burnside $\infty$-cate\-gory of finite $H$-sets -- is equivalent to the $G$-$\infty$-category of finitely generated free $G$-$E_{\infty}$-spaces. In other words, it is the free $G$-semiadditive $G$-$\infty$-category on one generator; that is, for any $G$-semiadditive $G$-$\infty$-category $A$, evaluation on the generator defines an equivalence of $G$-$\infty$-categories
\[\equivto{\underline{\Fun}_G^{\oplus}\left(\underline{A}^{\eff}(G),A\right)}{A}.\]
Here $\underline{\Fun}_G^{\oplus}$ is the $G$-$\infty$-category of $G$-coproduct-preserving functors.
\end{mainthm}

As suggested by work of Andrew Blumberg \cite{MR2286026}, the $G$-stability of a $G$-$\infty$-category can be defined as ordinary stability along with $G$-semiadditivity. Consequently, the two previous theorems, with some effort, together provide the following, also conjectured by Hill:

\begin{mainthm}\label{SpGuniv} The $G$-$\infty$-category $\underline{\Sp}^G$ of $G$-spectra -- whose value on an orbit $G/H$ is the $\infty$-category $\Sp^H$ of genuine $H$-spectra -- is the free $G$-stable $G$-$\infty$-category with all $G$-colimit on one generator; that is, for any $G$-stable $G$-$\infty$-category $E$, evaluation on the generator defines an equivalence of $G$-$\infty$-categories
\[\equivto{\underline{\Fun}_G^{L}\left(\underline{\Sp}^{G},A\right)}{A}.\]
\end{mainthm}

Going further, Hill also expected that the multiplicative norms of Hill--Hopkins--Ravenel would be part of a new type of structure -- a \emph{$G$-symmetric monoidal $G$-$\infty$-category.} In effect, a $G$-symmetric monoidal $G$-$\infty$-category is a $G$-$\infty$-category $C$ along with tensor product functors over finite $G$-sets. In particular, one has a functor
\[N^G\colon\fromto{C(G/e)}{C(G/G)},\]
which is exactly the desired multiplicative norm.

Work of Hill and Hopkins \cite{hillhopkins} has already laid out the idea of $G$-symmetric monoidal ordinary categories, but incorporating homotopy coherence into this sort of structure is a taller order. The situation is roughly analogous to the situation with the smash product in model categories of spectra: there are genuine obstructions to making a $G$-symmetric monoidal structure maximally compatible with a model category of genuine $G$-spectra. However, when we pass to the world of $\infty$-categories as in \cite{HA}, the situation becomes much cleaner: not only can one give an explicit, homotopy invariant construction of the smash product on the $\infty$-category $\Sp$ of spectra, but this smash product enjoys a universal property that characterizes it up to a contractible space of choices.

We bring exactly this kind of conceptual clarity (and technical power) to the study of homotopy coherent $G$-commutative structures in this text. We define the notions of \emph{$G$-$\infty$-operad} and \emph{$G$-symmetric monoidal $G$-$\infty$-category}. We find that $G$-products define $G$-symmetric monoidal structures on the $G$-$\infty$-category $\underline{\Cat}_{\infty,G}$ of $G$-$\infty$-categories and the $G$-$\infty$-category $\underline{\Top}_G$ of $G$-spaces. The $G$-commutative algebra objects of $\underline{\Cat}_{\infty,G}$ are precisely the $G$-symmetric monoidal $\infty$-categories, and the $G$-commutative algebra objects of $\underline{\Cat}_{\infty,G}$ are precisely the $G$-$E_{\infty}$-spaces.

Similarly, there is a $G$-subcategory $\underline{\Pr}^{L}_{G}\subset\underline{\Cat}_{\infty,G}$ of $G$-presentable $G$-$\infty$-cate\-gories and $G$-left adjoints. This too has a $G$-symmetric monoidal structure, but it is not given by $G$-products; rather, the $G$-commutative algebra objects of $\underline{\Pr}^{L}_{G}$ are precisely the $G$-symmetric monoidal $\infty$-categories that are presentable and in which the tensor product preserves $G$-colimits separately in each variable.

\begin{mainthm}\label{TopGunit} The $G$-$\infty$-category $\underline{\Top}_G$ is the unit in the $G$-symmetric monoidal $G$-$\infty$-category $\underline{\Pr}^{L}_{G}$. In particular, it admits an essentially unique $G$-symmetric monoidal structure in which the tensor product preserves $G$-colimits separately in each variable, which is given by the $G$-products.
\end{mainthm}

Even further, the full $G$-subcategory $\underline{\Pr}^{L}_{G,\textit{st}}\subset\underline{\Pr}^{L}_{G}$ spanned by the $G$-stable $G$-presentable $G$-$\infty$-categories inherits the $G$-symmetric monoidal structure, and we have the following.

\begin{mainthm}\label{SpGunit} The $G$-$\infty$-category $\underline{\Sp}^G$ is the unit object in the $G$-symmetric monoidal $G$-$\infty$-cate\-gory $\underline{\Pr}^{L}_{G,\textit{st}}$. In particular, it admits an essentially unique $G$-symmetric monoidal structure in which the tensor product preserves $G$-colimits separately in each variable.
\end{mainthm}

In particular, this provides a universal description of the Hill--Hopkins--Ravenel multiplicative norm. With some work, this even provides a universal characterization of an \emph{individual} norm functor.

\begin{mainthm}\label{normuniv} The norm functor $N^G\colon\fromto{\Sp}{\Sp^G}$ is the inital object of the $\infty$-category
\[\Fun^{\otimes}(\Sp,\Sp^G)\times_{\Fun^{\otimes}(\Top,\Sp^G)}\Fun^{\otimes}(\Top,\Sp^G)_{\Sigma^{G,\infty}_+\circ\Pi_G/},\]
where $\Pi_G\colon\fromto{\Top}{\Top^G}$ is the $G$-product, and $\Fun^{\otimes}$ denotes the $\infty$-category of symmetric monoidal functors.
\end{mainthm}

In this text, we completely realize Hill's vision, and we prove Theorems \ref{TopGuniv}--\ref{normuniv}.

\section*{Taking the $G$ out of Genuine} Formally, one may now note that the orbit category $\OO_G$ of the group $G$ plays a much more significant role in these results than does $G$ itself. In particular, although the $\infty$-category $\Sp^G$ can be obtained by taking the $\infty$-category $\Top_G$ and inverting the representation spheres, our viewpoint regards the role of representation spheres as incidental.

One is thus led to ask whether one might untether equivariant homotopy theory from dependence upon a group. (We thank Haynes Miller for the pun of the section heading.) That is, first, do Theorems \ref{TopGuniv}--\ref{normuniv} hold more generally? And, second, is there any value in proving them in greater generality? The answer to both questions turns out to be \emph{yes}.

Indeed, when one examines the proofs of the results above, one finds that the unstable results continue to hold when $\OO_G$ is replaced with any base $\infty$-category $T$. The stable results require only very mild conditions on $T$; in effect, one requires the analogue of the Mackey decomposition theorem in $T$ (``$T$ is \emph{orbital}'') and a condition that no nontrivial retracts exist (``$T$ is \emph{atomic}''). We can even extend this further, and define an \emph{incompleteness class} $R$ on the orbital $\infty$-category $T$; in effect, this serves to place limits on the classes of transfers that exist in the corresponding $\infty$-category of spectra.

As it happens, there are many examples that make this generality worthwhile. Here are a few.
\begin{enumerate}[1.]
\item As a mild extension of the example $\OO_G$, consider a family $\mathscr{F}$ of subgroups of $G$ such that if $K\leq L$ lie in $\mathscr{F}$, then any subgroup $H\leq G$ that is conjugate to a subgroup $H'$ such that $K\leq H'\leq L$ also lies in $\mathscr{F}$. Then the full subcategory $\OO_{G,\mathscr{F}}\subseteq\OO_G$ is also an atomic orbital category. Such categories (along with various inclusions of ``closed'' and ``open'' subcategories) appear naturally when we contemplate the isotropy separation sequence in equivariant stable homotopy theory.
\item Following Blumberg and Hill \cite{BH}, any incomplete $G$-universe $U$ gives rise to an incompleteness class $R_U$ on $\OO_G$, and this permits us to model $G$-spectra relative to $U$.
\item Furthermore, one can also work with orbit categories of profinite groups (where the stabilizers are required to be open) and locally finite groups (where the stabilizers are required to be finite). This provides extensions of equivariant stable and unstable homotopy theory to these contexts.
\item Any $\infty$-groupoid (= Kan complex) $X$ is atomic orbital. The corresponding $\infty$-category of \emph{$X$-spaces} is equivalent to the $\infty$-category of functors $\fromto{X}{\Top}$; likewise, the $\infty$-category of \emph{$X$-spectra} is equivalent to the $\infty$-category of functors $\fromto{X}{\Sp}$. In other words, $X$-spaces are local systems of spaces over $X$, and $X$-spectra are local systems of spectra over $X$. Consequently, this example actually recovers parametrized homotopy theory as studied by Peter May and Johann Sigurdsson \cite{MR2271789}; in fact, this example was the inspiration for our title.
\item Combining the previous example with the ur-example, for any $G$-space $X$, one can construct the \emph{total orbital $\infty$-category} \underline{X}. One sees that $\underline{X}$-spaces are local systems of $G$-spaces over $X$, and $\underline{X}$-spectra are local systems of $G$-spectra over $X$.
\item The cyclonic orbit $2$-category $\OO_{\cop}$, whose objects are $\QQ/\ZZ$-sets with finite stabilizers, whose $1$-morphisms are equivariant maps, and whose $2$-morphisms are certain intertwiners, is an orbital $\infty$-category \cite{cyclonic}. The corresponding homotopy theory of $B$-spectra is the homotopy theory of $S^1$-equivariant spectra relative to the family of finite subgroups. This is precisely the sort of equivariance that one sees on topological Hochschild homology \cite{tornado2}. To construct the homotopy category of \emph{cyclotomic} spectra, one forms the fixed points of this homotopy theory relative to the action of the monoid of \emph{open immersions} from $\OO_{\cop}$ into itself.
\item Generalizing the previous example are the \emph{multi-cyclonic orbit $2$-categories} which control torus-equivariance and \emph{multi-cyclotomic} structures, which appear naturally on higher forms of topological Hochschild homology \cite{multicyclonic}.
\item The $2$-category $\Gamma$ of finite connected groupoids and covering maps is atomic orbital. The corresponding homotopy theory of $\Gamma$-spectra is a variant of Stefan Schwede's global equivariant homotopy theory \cite{schwedeglobal}. To get \emph{exactly} Schwede's global equivariant homotopy theory (for finite groups) in our framework requires a larger orbital $\infty$-category of finite connected groupoids equipped with an incompletness class.
\item The category $\categ{Surj}_{\leq n}$ of finite sets of cardinality $\leq n$ and surjective maps is an atomic orbital category. This one is extremely strange, however, as it doesn't have much at all to do with any groups. Nevertheless, the third author shows in \cite{calcMack} that the corresponding homotopy theory of $\FF\SS_{\leq n}$-spectra is equivalent to the homotopy theory of $n$-excisive functors $\fromto{\Sp}{\Sp}$, generalizing Tom Goodwillie's classification of homogeneous functors. Indeed, it is the inclusion of $\categ{Surj}_{\leq n-1}$ into $\categ{Surj}_{\leq n}$, combined with the complmentary inclusion of $B\Sigma_n$ into $\categ{Surj}_{\leq n}$, that together reconstruct the recollement of $n$-excisive functors by $(n-1)$-excisive functors and $n$-homogeneous functors.
\item The $\infty$-categories $\categ{Surj}(\RR)_{\leq n}$ and $\categ{Surj}(\CC)_{\leq n}$ obtained from the topological categories of finite-dimensional inner product spaces (over $\RR$ and $\CC$, respectively) of dimension $\leq n$ and orthogonal projections are atomic orbital as well. Just as stable homotopy theory parametrized on the orbital categories $\categ{Surj}_{\leq n}$ ``controls'' the Goodwillie tower, so the stable homotopy theory parametrized on the orbital categories $\categ{Surj}(\RR)_{\leq n}$ ``controls'' the Weiss orthogonal calculus \cite{MR1321590}. Likewise, stable homotopy theory parametrized on the orbital categories $\categ{Surj}(\CC)_{\leq n}$ ``controls'' the unitary calculus. We hope to return to this point in future work.
\item Our framework also covers and extends a setting previously defined in work of Bill Dwyer and Dan Kan, Emanuel Dror Farjoun, and Boris Chorny and Bill Dwyer. In \cite{Farjoun}, Farjoun builds on work of \cite{DK} and defines a model structure on the category of diagrams of spaces indexed on a small category $I$, called the $I$-equivariant model structure, which depends on the ``$I$-orbits'': the diagrams $I\to \Top$ whose strict (= $1$-categorical) colimit is equal to a point. In particular if $I=G$ is a group these are precisely the $G$-orbits, and the resulting homotopy theory is the fixed-points model structure on $G$-spaces. Moreover Farjoun's construction admits an Elmendorf--McClure theorem, in the sense that the $I$-equivariant model structure is Quillen-equivalent to a presheaf category (on the orbit category when the orbits are either small or complete). This result was proved in different levels of generality in \cite{DK} and \cite{CD}, and in full generality in the more recent \cite{Chorny}. The category of $I$-orbits $\OO_I$ is an atomic orbital category, and by the above mentioned Elmendorf--McClure theorem, Farjoun's $I$-equivariant model category is equivalent to our homotopy theory of $\OO_I$-spaces. Our construction exhibits the $I$-equivariant homotopy theory as a fiber of a full-fledged $\OO_I$-category, thus enabling one to exploit the full theory of $\OO_I$-equivariant limits and colimits.
\end{enumerate}

Such a wealth of examples compels us to prove Theorems \ref{TopGuniv}--\ref{normuniv} in the generality atomic orbital $\infty$-categories, and, where possible, we develop elements of the theory in even greater generality.

\section*{Plan}

This text consists of nine Exposés:
\begin{enumerate}[I.]
\item We introduce the basic elements of the theory of \emph{parametrized $\infty$-categories} and functors between them. Following the lessons of \cite{BarShah}, these notions are defined as suitable fibrations of $\infty$-categories and functors between them. We give as many examples as we are able at this stage. Simple operations, such as the formation of opposites and the formation of functor $\infty$-categories, become slightly more involved in the parametrized setting, but we explain precisely how to perform these constructions. All of these constructions can be performed explicitly, without resorting to such acts of desperation as straightening. The key results of this Expos\'e are: (1) a universal characterization of the $T$-$\infty$-category of $T$-objects in any $\infty$-category, (2) the existence of an internal Hom for $T$-$\infty$-categories, and (3) a parametrized Yoneda lemma. \cite{Exp1}
\item We dive deep into the fundamentals of parametrized $\infty$-category theory in the second Expos\'e. In particular, we construct parametrized versions of join and slice, and use these to define \emph{parametrized colimits and limits} as well as \emph{parametrized left and right Kan extensions}. At the heart of this is the difficult but technically powerful result that, just as one may decompose colimits into coproducts and geometric realizations in $\infty$-category theory, similarly one may decompose parametrized colimits into parametrized coproducts and geometric realizations in the ordinary sense. This has the effect of elevating parametrized coproducts and products to a special status within the theory. Theorem \ref{TopGuniv} is proved (and generalized) here. \cite{Exp2}
\item We next introduce \emph{orbital $\infty$-categories}, along with a host of examples. There are actually different sorts of functor between orbital $\infty$-categories, and we taxonomize these according to certain algebro-geometric intuitions. For any orbital $\infty$-category $T$, we have a corresponding $\infty$-category $\Sp^{T}$ (even $T$-$\infty$-category) of \emph{$T$-spectra}, which under our algebro-geometric analogy corresponds roughly to an $\infty$-category of ``quasicoherent sheaves on $T$.'' The different sorts of functors between orbital $\infty$-categories induce suitable functors between the $\infty$-categories of $T$-spectra, and these behave as the names suggest. Perhaps most importantly, \emph{closed immersions} of orbital $\infty$-categories admit \emph{open complements}, and these two functors induce a \emph{recollement} of the corresponding $\infty$-category of spectra; this is how one obtains the \emph{isotropy separation sequence} and generalizations thereof. With a little care, we are able to extend all this to the context of an orbital $\infty$-category equipped with an incompleteness class. \cite{Exp3}
\item In the fourth Expos\'e, we define \emph{semiadditive parametrized $\infty$-categories}, and we prove Theorem \ref{Aeffuniv}. Then we use the work of Expos\'e II to show that parametrized stability can be expressed as ordinary stability combined with parametrized semiadditivity. This now makes it possible to prove Theorem \ref{SpGuniv}. \cite{Exp4}
\item Next, we introduce the notion of \emph{parametrized Waldhausen $\infty$-categories}. We show that the algebraic $K$-theory of a Waldhausen $T$-$\infty$-category naturally carries the structure of a $T$-spectrum. \cite{Exp5}
\item From here, we move toward the algebraic structures in parametrized higher category theory. We introduce the notions of \emph{$T$-$\infty$-operad} and \emph{$T$-symmetric monoidal $\infty$-category} for an orbital $\infty$-category $T$, and we offer up numerous examples. Perhaps most importantly, parametrized $\infty$-categories with all $T$-coproducts (or, dually, $T$-products) inherit canonical $T$-symmetric monoidal structures. \cite{Exp6}
\item In the seventh Expos\'e, we prove that when $T$ is an atomic orbital $\infty$-category, the $T$-$\infty$-category $\underline{\Pr}_T^L$ of $T$-presentable $T$-$\infty$-categories admits a $T$-symmetric monoidal structure analogous to the symmetric monoidal structure on presentable $\infty$-categories. Theorem \ref{TopGuniv} then implies, more or less directly, Theorem \ref{TopGunit}. Moreover, $T$-stable $T$-presentable $T$-$\infty$-categories form a symmetric monoidal localization of $\underline{\Pr}_T^L$, and the localization is given by tensoring with $\underline{\Sp}^T$. Theorem \ref{SpGunit} follows immediately. Moreover, one deduces a different universal property of $\underline{\Sp}^T$, which is that it is, in effect, the result of inverting the analogues of the permutation representation spheres in the $T$-symmetric monoidal $T$-$\infty$-category $\underline{\Top}_G$. From this, we are able to deduce the universal property of the norm (Theorem \ref{normuniv}). \cite{Exp7}
\item In the penultimate Expos\'e, we introduce the $T$-$\infty$-category $\underline{\Mod}(A)$ of modules over a $T$-$E_{\infty}$-algebra $A$ (for an atomic orbital $\infty$-category $T$). We show that it is $T$-symmetric monoidal, and we describe how it transforms in both $A$ and $T$. \cite{Exp8}
\item Finally, we return to the subject of equivariant algebraic $K$-theory, where we show that the equivariant algebraic $K$-theory of a $T$-symmetric monoidal Waldhausen $T$-$\infty$-cate\-gory admits the natural structure of a $T$-$E_{\infty}$ ring spectrum. This applies not only in the field case of the beginning of this introduction, but also to those forms of equivariant algebraic $K$-theory that arise in the work of Dustin Clausen, Akhil Mathew, Niko Naumann, and Justin Noel \cite{CMNN} as well as the nascent subject of equivariant (derived) algebraic geometry. \cite{Exp9}
\end{enumerate}

\section*{Acknowledgments}

This text answers a web of questions and conjectures posed by Mike Hill and Mike Hopkins. We thank them for presenting us with these questions in the afterglow of their solution with Doug Ravenel to the Kervaire Invariant Problem, as well as for their encouragement as we developed this project. We thank Mike Hill in particular for his vision in presenting his extensive program to build up the foundations of equivariant derived algebraic geometry.

Additionally, we thank the other participants of the Bourbon Seminar at MIT -- Lukas Brantner, Peter Haine, Marc Hoyois, and Akhil Mathew, for the many, many conversations we had with them over the course of the completion of this text. We also thank John Berman for his very surprising insights.


\bibliographystyle{amsplain}
\bibliography{Gcats}

\end{document}